\documentclass[twocolumn]{autart}    

\usepackage[dvipsnames]{xcolor}
\usepackage{textcomp}
\usepackage{comment}
\usepackage{graphicx}
\usepackage{amssymb}
\usepackage{gensymb}
\usepackage{amsmath}

\usepackage[labelfont={bf},font=small]{caption}
\usepackage[none]{hyphenat}
\usepackage{mathtools, cuted}
\usepackage[noadjust, nobreak]{cite}
\usepackage{tabularx}
\usepackage{amsmath}
\usepackage{float}
\usepackage{pifont}
\usepackage[]{placeins}
\usepackage{makecell}
\usepackage{tikz}
\usepackage[framemethod=tikz]{mdframed}
\usepackage{afterpage}
\usepackage{stfloats}
\usepackage{atbegshi}
\AtBeginShipout{\handlethispage}
\usepackage{acronym}
\usepackage{hyperref}

\usepackage[round]{natbib}

\usepackage{amsfonts}
\usepackage{eso-pic}

\usepackage[font=small,skip=0pt]{caption}
\usepackage[acronym]{glossaries}

\usepackage{algpseudocode}

\usepackage[mathscr]{eucal}

\newcolumntype{Y}{>{\centering\arraybackslash}X}

\newcommand{\handlethispage}{}
\newcommand{\discardpagesfromhere}{\let\handlethispage\AtBeginShipoutDiscard}
\newcommand{\keeppagesfromhere}{\let\handlethispage\relax}
\newcommand{\RNum}[1]{\uppercase\expandafter{\romannumeral #1\relax}}

\hypersetup{
colorlinks=true, 
breaklinks=true, 
urlcolor= black, 
citecolor=TealBlue,
linkcolor= TealBlue, 
bookmarksopen=TealBlue,
pdftitle={Learning-based model predictive control with moving horizon state estimation for autonomous racing}, 
pdfauthor={Y. Kebbati},
}

\edef\endfrontmatter{%
  \unexpanded\expandafter{\endfrontmatter}
  \noexpand\endNoHyper 
}

\usepackage{background}

\SetBgAngle{-90}
\SetBgScale{8}
\SetBgColor{gray!60}
\SetBgContents{%
  \begin{tikzpicture}
    \matrix [column sep =0.4cm,ampersand replacement=\&]
      { \node [rotate=90] {\hspace{2.5cm} Preprint}; \& \node [rotate=90] {};\\};
  \end{tikzpicture}}

\begin{document}

\begin{frontmatter}

\title{Learning-based model predictive control with moving horizon state estimation for autonomous racing} %

\thanks{Corresponding author Yassine Kebbati.}

\author[1]{Yassine Kebbati}\ead{yassine.kebbati@univ-evry.fr},    
\author[3]{Andreas Rauh}\ead{andreas.rauh@uni-oldenburg.de}, 
\author[1]{Naima Ait-Oufroukh}\ead{naima.aitoufroukh@univ-evry.fr},               
\author[1]{Dalil Ichalal}\ead{dalil.ichalal@univ-evry.fr},
\author[1,2]{Vincent Vigneron}\ead{vincent.vigneron@univ-evry.fr}  

\address[1]{IBISC-EA4526, Paris-Saclay University, Evry, France}  
\address[2]{School of Applied Sciences (FCA),UNICAMP, Limeira, Brazil}    
\address[3]{Distributed Control in Interconnected Systems, Carl von Ossietzky University, Oldenburg, Germany}


\begin{abstract}
This paper addresses autonomous racing by introducing a real-time nonlinear model predictive controller (NMPC) coupled with a moving horizon estimator (MHE). The racing problem is solved by an NMPC-based off-line trajectory planner that computes the best trajectory while considering the physical limits of the vehicle and circuit constraints. The developed controller is further enhanced with a learning extension based on Gaussian process regression that improves model predictions. The proposed control, estimation, and planning schemes are evaluated on two different race tracks.
\end{abstract}

\end{frontmatter}

\section{Introduction} \label{sec:introduction}
Self-driving vehicles have seen tangible advancements in recent years, especially in the fields of decision-making, planning, and control. Many advanced driver assistance systems (ADAS), such as cruise control, emergency braking, and lane-keeping, have been deployed in modern commercial vehicles for quite some time. However, challenging driving situations call for more sophisticated control algorithms. Thus, several new approaches regarding vehicle dynamics control have recently been introduced in the literature, ranging from classic model-based to model-free and data-driven control strategies.   

For instance, \cite{Corno2020} addressed vehicle longitudinal control with a linear parameter varying (LPV) $H_{\infty}$ strategy for fast driving under evasive maneuvers. Similarly, \cite{xu2018accurate} integrated road slope, speed profile, and vehicle longitudinal dynamics in an augmented optimal preview controller for speed regulation. A self-adaptive PID with genetic algorithms and neural networks has been proposed in \cite{kebbati2021optimized} for adaptive speed regulation. \cite{rokonuzzaman2021customisable} introduced a data-driven MPC longitudinal controller that learns from human driving demonstrations. 

On the other hand, model-free techniques for lateral control have been studied in \cite{bojarski2016end}, \cite{7963716} and \cite{Jhung2018}. The authors used camera images to ensure end-to-end steering control based on closed-loop feedback with a convolutional neural network (CNN). In a similar approach, \cite{Eraqi2017} combined recurrent neural networks (RNN) with CNN to predict steering controls. However, the application of such techniques remains a challenge due to safety requirements and because these methods lack interpretability and cannot guarantee stability. Alternatively, model-free techniques can be used to enhance model-based control strategies. For instance, an adaptive MPC for path tracking was proposed by \cite{kebbati2021optimized2}, which was further improved in \cite{kebbati2021neural} by including neural networks and adaptive neuro-fuzzy inference systems to the adaptation approach. \cite{salt2021autonomous} addressed lane-keeping by an LPV-MPC controller with a dual-rate extended Kalman filter for state estimation. A recent review of the most widely used techniques for lateral control is provided in \cite{kebbati2022lateral}.

Although the above-mentioned strategies handle both lateral and longitudinal control tasks separately, full driving autonomy requires simultaneous lateral and longitudinal vehicle control. In this regard, \cite{attia2014combined} and \cite{Chebly2019} proposed a combined lateral and longitudinal control strategy in which both control tasks are coordinated to achieve full autonomy. In a similar way, \cite{kebbati2022coordinated} introduced a coordinated control strategy using LPV-MPC for lateral control with a Particle Swarm Optimized Proportional Integral Derivative controller (PSO-PID) for speed regulation.

The aforementioned autonomous driving techniques are meant for different driving scenarios, the driving task can range from simple highway or urban driving to more complex scenarios such as autonomous racing. In such applications, the objective is to go as fast as possible, where the vehicle racing autonomously has to operate at its handling limits leading to a more challenging control problem. Knowing the physical limits of the vehicle is important for such driving scenarios, which implies the necessity of an accurate model. Autonomous racing, as a sub-field of autonomous driving, has recently attracted considerable research attention \cite{Liniger2015,verschueren2016time,Kabzan2019a,tuatulea2020design}, and \cite{Alcala2020}. In this regard, the literature favors MPC over the other model-based control techniques, especially for lateral control applications. MPC is known for its good performance and ability to handle autonomous driving since it allows to accurately predict the required driving behaviour even in challenging situations. Furthermore, MPC is known for its unique features such as imposing safety and physical constraints and handling multi-input multi-output (MIMO) systems, in addition to its predictive nature that can foresee the vehicle's behaviour beforehand. However, the latter is challenging due to the compromise between precision and fast computation \cite{Alcala2020,verschueren2014towards}. On the one hand, using complex nonlinear models yields accurate results at the expense of slow computations, which is unsuitable for real-time applications. On the other hand, simpler models like linear dynamic or kinematic bicycle models run relatively fast but lack precision.

In addition to control, race line planning is a very important task in autonomous racing. Several techniques exist for this purpose, sampling-based planning algorithms tend to quickly find feasible trajectories, but they may not be optimal. These are divided into graph-based and incremental search algorithms. The best-known graph-based methods are the Dijkstra and A* algorithms, which use heuristic information to search for cost-optimal paths within a discretized configuration space of the vehicle, known as the graph \cite{karur2021survey}. On the other hand, the rapidly exploring random trees (RRT) and RRT* are the most widely used incremental search planning algorithms \cite{heilmeier2019minimum}. RRT has been investigated by \cite{hwan2013optimal} to generate a time-optimal trajectory for a 180$\degree$ curve. Alternatively, optimization-based methods such as optimal control can find the fastest race line while taking into account the physical limits of the vehicle. \cite{verschueren2016time}, and \cite{Liniger2015} proposed a one-level approach combining control and planning in a single optimization problem. \cite{Alcala2020} proposed using a two-level approach to separate control and planning, claiming that one-level methods complicate the control problem and may require more computation resources. 
In this paper, we propose a complete setup for autonomous racing, where the control part is ensured using a nonlinear model predictive controller (NMPC) with a low-order kinematic model that includes longitudinal dynamics. Furthermore, a moving horizon estimator (MHE) is developed for real-time state estimation to complement the controller. The controller seeks to track the optimal race line of the circuit, which is planned offline using an NMPC with a high-fidelity nonlinear vehicle dynamics model. The main contribution of this paper is a learning-based MPC control algorithm that uses a low-order kinematic model to significantly reduce computation efforts while guaranteeing accurate tracking. A learning extension that uses Gaussian process ($\mathscr{GP}$) regression is proposed to correct the low-order model mismatch and enhance control performance. 

This article is divided as follows: Section 2 addresses the modeling of the vehicle dynamics for control, planning, and estimation. Section 3 details the design of the NMPC controller, the MHE state estimator, the race line planner, and the learning extension of the controller. Evaluation results of the estimation, control performance, and learning approach are presented and analyzed in Section 4. Finally, Section 5 summarizes the paper with conclusions and gives directions for future work. 

\section{Vehicle Modeling}
\begin{figure}[tb]
\centering
\includegraphics[width=7.75cm,height=5.5cm]{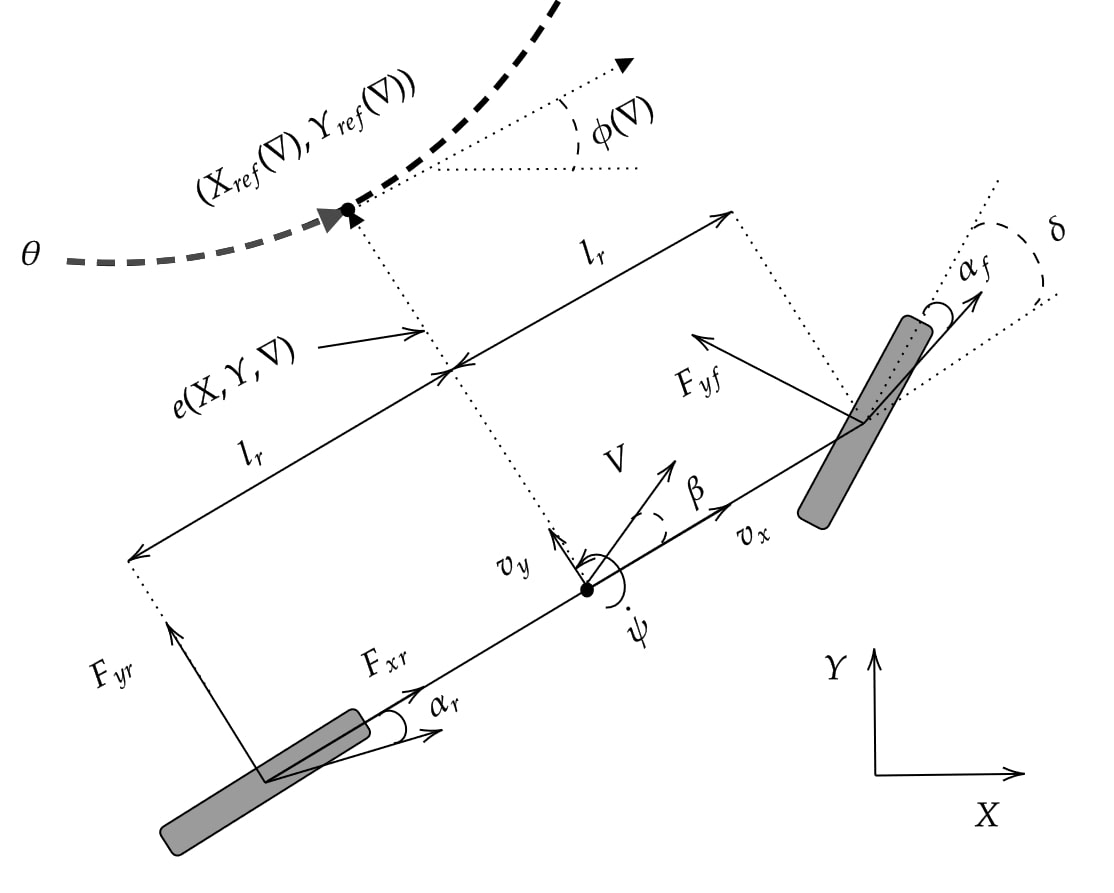}
\caption{Bicycle model.}
\label{fig1}
\end{figure}

The vehicle in question is a $1/10$ scale electric race car modeled as a rigid body with mass $m$ and inertia $I_z$. A single-track bicycle model can be formulated assuming the car is symmetrical and considering in-plane motions only (see Fig. \ref{fig1}). This modeling approach is one of the most common ways to model racing cars \cite{Liniger2015,Kabzan2019a,tuatulea2020design}, where the state vector is defined as $\breve{x} = [X,Y,\psi,v_x,v_y,\omega]^T$. Parameters $X$ and $Y$ are the global coordinates of the car's position, and $\psi$ represents its orientation. The linear lateral and longitudinal velocities of the car are represented by $v_x$ and $v_y$, respectively, and $\omega$ is the car's angular velocity. Furthermore, the control input is defined as $u = [\delta, \ D]^T$, ith $\delta$ being the front steering angle and $D$ the motor's duty cycle. The model is governed by the following ordinary differential equations (ODEs):

\begin{subequations}
\label{eq:0}
\begin{align}
&\Dot{X} = v_x \cos(\psi) - v_y \sin(\psi),\label{eq01}\\
&\Dot{Y} = v_x \sin(\psi) + v_y \cos(\psi),\label{eq02}\\
&\Dot{\psi} = \omega, \label{eq03}\\
&\Dot{v}_x = \frac{F_{x,r} - F_{y,f} \sin(\delta) +mv_y\omega}{m},\label{eq04}\\
&\Dot{v}_y = \frac{F_{y,r} - F_{y,f} \cos(\delta) -mv_x\omega}{m},\label{eq05}\\
&\Dot{\omega} = \frac{F_{y,f}l_f\cos(\delta) - F_{y,r}l_r}{I_z}.\label{eq06}
\end{align}
\end{subequations}
The tire-road interactions are modeled using the simplified Pacejka tire model, which is realistic enough for high-speed maneuvers at vehicle handling limits \cite{Alcala2020}.

\begin{subequations}
\label{eq:9}
\begin{align}
&F_{y,f} = D_f \sin(C_f \arctan(B_f \alpha_f)),\label{eq91}\\
&F_{y,r} = D_r \sin(C_r \arctan(B_r \alpha_r)),\label{eq92}\\
&F_{x,r} = (C_{m1} - C_{m2} v_x)D - C_{ro} - C_d v_x^2\label{eq922}\\
&\alpha_f = \delta -\arctan\left( \frac{v_y + l_f \Dot{\psi}}{v_x} \right),\label{eq93}\\
&\alpha_r = \arctan\left( \frac{l_r \Dot{\psi} - v_y}{v_x} \right).\label{eq94}
\end{align}
\end{subequations}
Equations (\ref{eq91}-\ref{eq92}) define the tire lateral forces. The parameters $B$, $C$ and $D$ determine Pacejka's semi-empirical curve \cite{Pacejka2008}, which fits empirical data to model the tire dynamics. The subscripts $(_{f,r})$ refer to the front and rear wheels, respectively. The term $\alpha$ denotes the side-slip angle of the corresponding (front/rear) wheel. Parameter $l$ is the distance between the vehicle's center of gravity (CG) and the corresponding wheel's axle. Finally, equation (\ref{eq922}) represents the rear traction force, which is modeled by a simplified representation of a DC motor with air drag and rolling resistance. The coefficients for rolling resistance and drag are represented by $C_{ro}$ and $C_d$, and the terms $C_{m_{(1,2)}}$ are motor specific coefficients.

The model described above is very accurate with high fidelity but expensive in terms of computation for real-time applications, most works report a minimum average computation time of around 50ms \cite{Kabzan2019a,Alcala2020}. Hence, this model is used in this paper for simulating vehicle behaviour and validating the controller. At this stage and to simplify control design under a low computational burden, we make use of a low-order vehicle model that efficiently runs within real-time optimization. For this model, load changes, pitch, roll dynamics, and road tire interactions are neglected, resulting in a slip-free model \cite{verschueren2014towards}, summarized in the following equations:

\begin{subequations}
\label{eq:1}
\begin{align}
&\Dot{X} = v \cos(\psi + g_1 \delta),\label{eq11}\\
&\Dot{Y} = v \sin(\psi + g_1 \delta),\label{eq12}\\
&\Dot{\psi} = v \delta g_2,\label{eq13}\\
&\Dot{v} = (C_{m_1}-C_{m_2}v)D - C_{r_2} v^2 - C_{r_1} -(v \delta)^2 g_1^2 g_2.\label{eq14}
\end{align}
\end{subequations}
where $g_1$ and $g_2$ are geometric parameters properly identified along with the other parameters in table \ref{tab1}. The velocity of the car is approximated, assuming small steering angle $\delta$, as follows  \cite{verschueren2014towards}:
$$v_x \approx v, \ \ v_y \approx vg_1\delta.$$
Thus, the vehicle's velocity reads: $v^2 = v_x^2 + v_y^2$.
Observing model (\ref{eq:1}), one can notice that equations (\ref{eq11}-\ref{eq13}) define the vehicle's kinematics, while equation (\ref{eq14}) represents the longitudinal dynamics that depend on the duty cycle $D$ and the resistive force.

\begin{table}[tb]
\caption{Modeling parameters.} 
\label{tab1}
\centering
\begin{tabular}{c c c c c  }

\hline
Parameter & Unit & Value & Parameter  & Value \\ [0.5ex] 

$m$ & $\mathrm{kg}$  & $1.98$ & $B_f$   & $29.5$\\
 
$I_z$ & $\mathrm{kg.m^2}$  & $0.1217$ & $C_f$   & $0.087$\\

$l_f$ & $\mathrm{m}$  & $0.125$ & $D_f$  & $42.53$ \\

$l_r$ & $\mathrm{m}$  & $0.125$ & $B_r$  & $26.97$ \\

$C_{m_1}$ & $\mathrm{m/s^2}$  & $12$ & $C_r$  & $0.163$\\

$C_{m_2}$ & $\mathrm{s^{-1}}$  & $2.17$ & $D_f$  & $161.59$\\

$C_{r_1}$ & $\mathrm{m^{-1}}$  & $0.6$ & $g_1$ & $(l_r/(l_r+l_f))$ \\

$C_{r_2}$ & $\mathrm{m/s^2}$  & $0.1$ & $g_2$  & $(1/(l_r+l_f))$ \\[1ex] 

\end{tabular}
\end{table}

\section{Controller Design}

\subsection{Model Predictive Control}\label{nmpc}
To ensure accurate trajectory tracking, an NMPC controller is formulated to handle both steering and acceleration controls. The NMPC controller uses vehicle model (\ref{eq:1}), discretized by the sampling time $T_s$, to predict the vehicle's behaviour over a certain prediction horizon $N$. Then, it generates a sequence of optimal control actions, that minimize tracking errors by solving a constrained nonlinear optimization problem (NLP). Because of the receding horizon principle, only the first term of the optimal control actions is applied, and the whole process is repeated at each iteration. In this regard, the NLP problem reads as follows:

\begin{subequations}
\label{eq:2}
\begin{align}
\label{eq21}
\min_{\Delta U_k} \sum_{k=1}^{N}& \big\{||r_{k} - x_{k}||^2_Q  + || \Delta u_{k}||^2_R \big\} + ||x_{k+N+1}||_P^2. \\
\label{eq22}
s.t:\ & x_{k+1} = f(x_k, u_k), \\
\label{eq23}
    & u_k = u_{k-1} + \Delta u_k,\\
\label{eq24}
    &  u_k \in U, \\
\label{eq25}
    &  \Delta u_k \in \Delta U,\\
\label{eq26}
    &  x_k \in  X,\\
\label{eq27}
    &  x_{k_{\{1,2\}}} \leq  B_R(k),\\
\label{eq28}
    &  B_L(k) \leq  x_{k_{\{1,2\}}}.
\end{align}
\end{subequations}
The states and control inputs are defined as $x = [X, Y, \psi, v]^T$, $u = [\delta, D]^T$, and $N$ is the previously defined prediction horizon. The positive semi-definite diagonal matrices $Q$ and $R$ are weighting matrices that penalize the states and the control effort, $r_k$ is the reference trajectory with global coordinates $X$ and $Y$, and the last term of (\ref{eq21}) is the terminal cost used to increase stability. Generally speaking, the asymptotic stability of MPC with a quadratic stage cost is ensured by adding a terminal cost and a terminal set. Otherwise, a sufficiently long prediction horizon is needed invoking cost controllability to guarantee asymptotic stability as illustrated by \cite{russwurm2021mpc}, who proposed using the theory of barriers with MPC cost controllability to guarantee MPC asymptotic stability. \cite{muller2017quadratic} further argues that MPC stability can be satisfied if the system satisfies a local controllability condition or in case its linearization at the origin is stabilizable. In addition, \cite{mayne2000constrained} states that adding a terminal cost to the MPC formulation and using a sufficiently long prediction horizon ensures MPC stability without the need for terminal constraints. The system dynamics are defined in (\ref{eq22}), while (\ref{eq24}-\ref{eq26}) represent the constraints on the control inputs, their rates, and the states respectively. To keep the vehicle inside the circuit, hard constraints (\ref{eq27},\ref{eq28}) are imposed on the vehicle's position ($X$ and $Y$ coordinates), where the terms $B_L$ and $B_R$ are $3^{rd}$ order spline approximations of the circuit's boundaries. 

Autonomous racing is a complicated task in which the vehicle must operate at the limit of its handling capabilities. Moreover, the trajectory to be followed needs to be a time-optimal race line that allows the vehicle to go as fast as possible while operating at its dynamic limits, meaning that the vehicle is exposed to higher lateral and longitudinal accelerations compared to normal driving. The planning of such a race line is explained in the following subsection.   

\subsection{Race Line Planning}
The objective of path planning for racing is to generate an optimal race line that minimizes lap time, this time optimality can be seen as the maximum feasible progress along the racing track. The latter is achieved by exploiting spatial coordinates instead of time domain coordinates, meaning that vehicle properties such as positions, angles, and so on are represented in the spatial domain. To generate the race line, an offline NMPC with the previously introduced full vehicle model (\ref{eq:0}) is formulated in the spatial domain where the center line of the circuit is used as a measure of progress. The key point here is to use a long prediction horizon and low weights on the tracking error to maximize progress at the expense of tracking accuracy. Using the full vehicle dynamic model with a long prediction horizon requires very long computations, therefore, this phase is computed offline to generate the race line. 
The racing path is parameterized through offline-fitted polynomial interpolation by the arc length $\theta \in [0,\ L]$, with $L$ being the total length. Hence, any point of the parameterized center-line $(X_{ref}(\theta),Y_{ref}(\theta))$ (see Fig. \ref{fig1}) can be obtained using the aforementioned parameterization, that is by evaluating the fitted polynomial for the argument $\theta$. At this point, one can minimize the spatial tracking error defined as the deviation of the car's position from the reference point $(X_{ref}(\theta), Y_{ref}(\theta))$ by using a projection operator as follows \cite{Liniger2015}:

\begin{equation}
\label{eq:3}
    \nabla(X,Y) \triangleq \min_{\theta}\left[(X-X_{ref}(\theta))^2+(Y-Y_{ref}(\theta))^2\right].
\end{equation}
Using the definition $\nabla  \triangleq \nabla(X,Y)$, the orthogonal distance between the car and the reference path is given by the following expression:

\begin{equation}
\label{eq:4}
\begin{split}
        e(X,Y,\nabla ) &\triangleq \ \sin(\phi(\nabla ))(X-X_{ref}(\nabla ))\\
    &-\cos(\phi(\nabla ))(Y-Y_{ref}(\nabla )).
\end{split}
\end{equation}
where:

\begin{equation}
\label{eq:5}
    \phi(.) \triangleq \arctan \left\{ \frac{\partial Y_{ref}(.)}{\partial X_{ref}(.)}\right\}.
\end{equation}

The objective now becomes a matter of minimizing the orthogonal distance $e(X, Y, \nabla)$ while maximizing the racing progress represented by the projection $\nabla$ over a long prediction horizon $N_p$. Hence, the optimization problem can be posed as follows:

\begin{subequations}
\label{eq:7}
\begin{align}
\label{eq71}
\min_{u,\theta} \sum_{k=1}^{N_p} & \big\{ ||e(X_k,Y_k,\nabla)||_q^2 \big\} - r\nabla_{N_p}.\\
\label{eq72}
s.t:\ & x_{k+1} = g(x_k, u_k), \\
\label{eq73}
    &  u_k \in U, \\
\label{eq74}
    &  x_k \in  X,\\
\label{eq75}
    & B_L(k) \leq x_{k_{\{1,2\}}} \leq  B_R(k).
\end{align}
\end{subequations}

The trade-off between progress and reference tracking is obtained through the weights $q$ and $r$. Equation (\ref{eq72}) is the discrete version of the full vehicle dynamics given in (\ref{eq:0}) and the rest is the same as in problem (\ref{eq:2}). Notice that the objective function in (\ref{eq71}) depends on the projection operator (\ref{eq:3}) and is not real-time feasible. Hence, this part is computed offline to generate the race line which will be used in the real-time control problem (\ref{eq:2}) with the simplified model (\ref{eq:1}).    

\subsection{Moving Horizon Estimation }
\begin{figure}[htb]
\centering
\includegraphics[width=7.75cm,height=5.5cm]{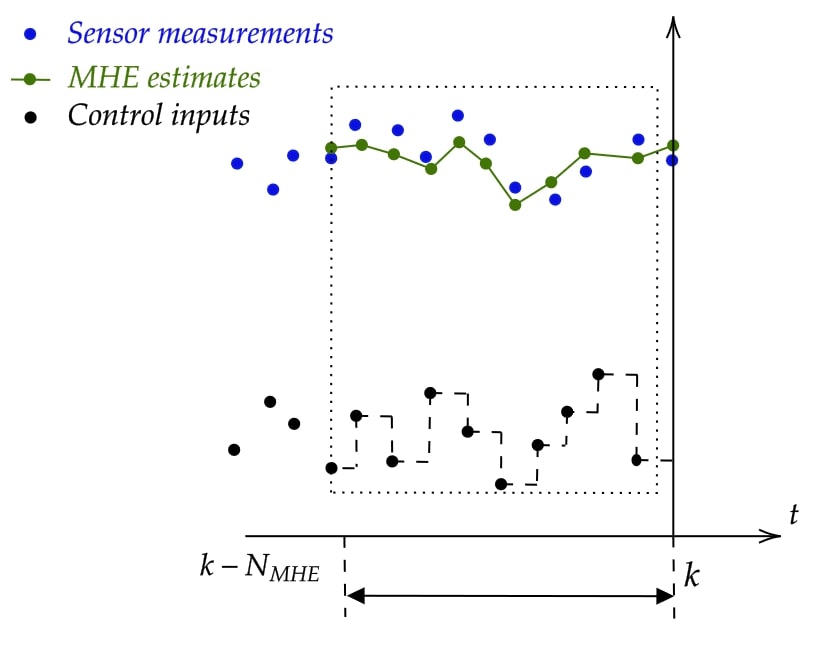}
\caption{Moving horizon estimation process.}
\label{fig2}
\end{figure}

Autonomous vehicles are equipped with a variety of sensors that can measure most of the dynamic states of the vehicle, but information coming from cheap sensors is always noisy and contains inaccuracies which calls for rigorous post-processing. Hence, MHE becomes a very handy solution since it makes use of past measurements to produce accurate estimates of the states. The goal of the moving horizon estimator (MHE) is to predict the dynamic states of the vehicle for future control iterations. This is achieved by solving a constrained optimization problem online over a fixed estimation horizon. Although the Kalman filter is a mature and well-known state estimation technique, it usually is limited to linear unconstrained systems and cannot handle nonlinear constrained problems. The obvious advantage of MHE is its ability to optimally handle constrained nonlinear systems using a finite horizon of measurement data, which allows it to bound the optimization problem \cite{haseltine2005critical}. The MHE can be seen as the dual problem of the NMPC controller discussed in section \ref{nmpc}, where at each iteration $k$, the MHE acquires past measurements over the estimation horizon $N_e$ and then finds the states' trajectory that best fits the considered window of measurements (as illustrated in Fig. \ref{fig2}), which is achieved by solving the following NLP problem: 
\begin{subequations}
\label{eq:8}
\begin{align}
\label{eq81}
\min_{x,u}\sum_{i=k-N_e}^k&||\tilde{y}(i)-h(x(i))||^2_V + \sum_{i=k-N_e}^{k-1}||\tilde{u}(i)-u(i)||^2_W.\\ 
\label{eq82}
s.t:\ & x(i+1) = f(x(i), u(i) + \omega_u(i)),  \\
\label{eq83}
    & y(i) = h(x(i)) + \omega_y(i),  \\
\label{eq84}
    &  u_k \in U, \\
\label{eq85}
    &  x_k \in  X.
\end{align}
\end{subequations}
with equation (\ref{eq82}) being the state propagation model (same system dynamics introduced in (\ref{eq:1})), and (\ref{eq83}) representing the sensor measurement model. Here $\omega_y \thicksim \mathcal{N}(0,\sigma_y) $  and $\omega_u \thicksim \mathcal{N}(0,\sigma_u)$ represent the output and input measurement noise with their respective covariance matrices $\sigma_y$ and $\sigma_u$. The terms $\Tilde{y}$ and $\Tilde{u}$ are the measured system outputs and inputs, respectively, and $||.||_{\{V,W\}}$ represents the Euclidean norm where $V$ and $W$ are positive-definite weighting matrices.

\subsection{Learning Extension}

The learning extension aims to improve MPC model predictions to enable the vehicle to drive more aggressively and increase its racing capability. Hence, the aim is to try to find the mismatch between the high fidelity model (\ref{eq:0}) that represents the vehicle, and the low order model (\ref{eq:1}) used in the MPC design. The model mismatch can be expressed as the error between MPC predicted states and actual vehicle states coming from the high-fidelity model. Let us define $F_\text{vehicle}$ as the high-fidelity model and $F_\text{mpc}$ as the low-order model used for MPC predictions. The model mismatch accounts for the four common states between the two models and has the following form: $[e_X, e_Y, e_{\psi}, e_v, 0, 0]^T$. Therefore, a Gaussian process regressor is learned for each state using recorded data from previous racing laps performed by the regular NMPC controller.

Gaussian processes are a collection of random variables having a joint Gaussian distribution. Consider noisy observations coming from a model function $f : \mathbb{R}^n \rightarrow \mathbb{R}$ through a Gaussian noise model $o = f(x) + \mathcal{N}(0,\sigma_n^2)$ with $x \in \mathbb{R}^n$. The Gaussian process of $o$ is defined by the mean function $\mu(x)$ and the covariance function $k(x,x')$ \cite{williams2006gaussian}. 
\begin{equation}
\label{eq:11}
\begin{split}
    \mu(x,\zeta) & = \mathbb{E}[f(x)],\\
    k(x,x',\zeta) & = \mathbb{E}[(f(x) - \mu(x))(f(x')-\\
    &\mu(x'))] + \sigma_n^2 \delta(x,x'). 
\end{split}
\end{equation}
where $\delta(x,x')$ is known as the Kronecker delta function, and $\zeta$ is a hyperparameter vector that parameterizes the mean and covariance functions.

\begin{figure}[htb]
\centering
\includegraphics[width=7.75cm,height=5.5cm]{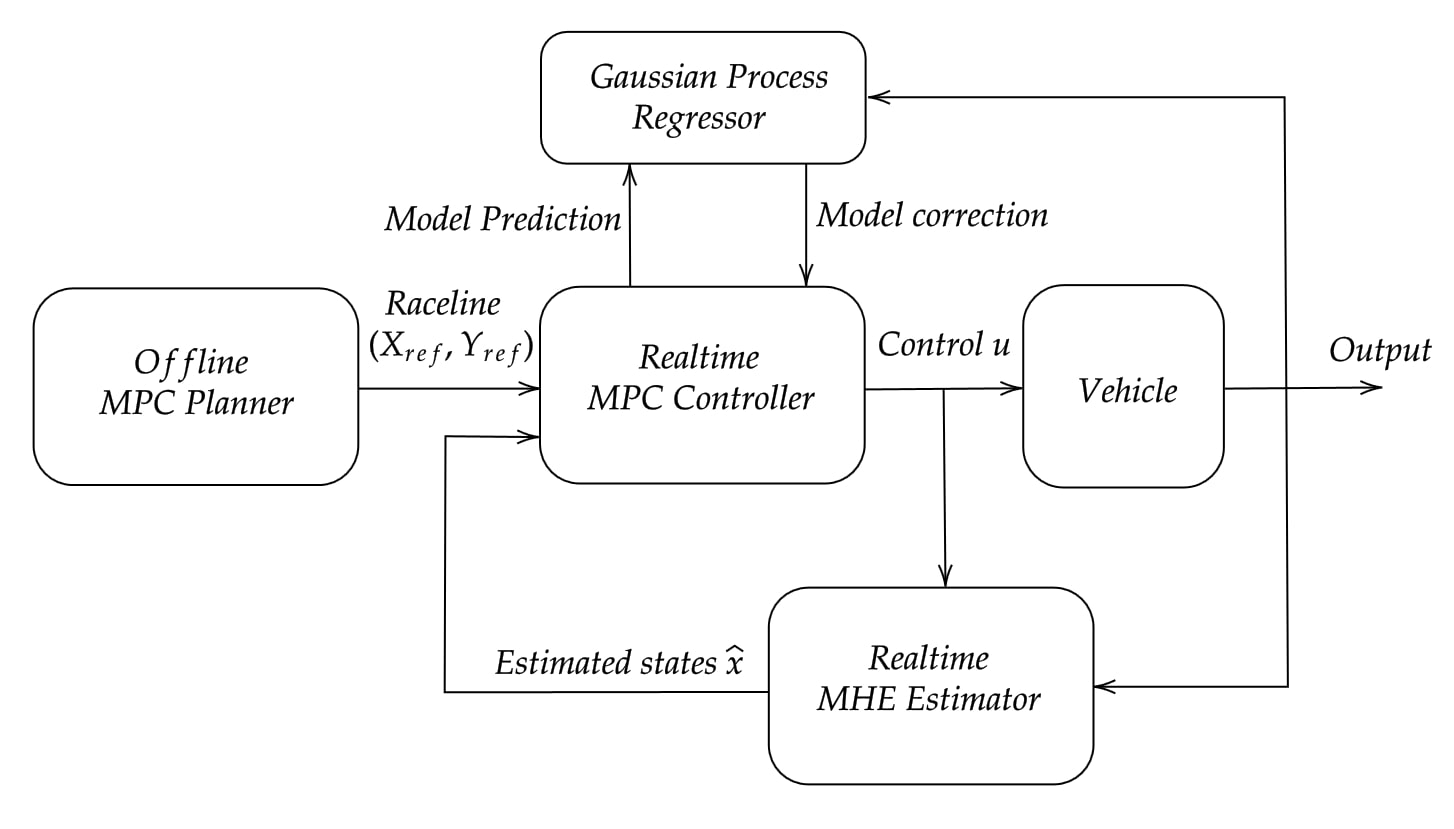}
\caption{Overall control approach.}
\label{fig3}
\end{figure}

The training data are obtained from both $F_\text{vehicle}$ and $F_\text{mpc}$ in addition to the corresponding control inputs. The model mismatch can be expressed as:
\begin{equation}
\label{eq:12}
F_\text{error}(x_k,u_k) = F_\text{vehicle\{1:4\}}(x_{k+1}) - F_\text{mpc}(x_k,u_k).
\end{equation}
where $k\in \{0,1,...,T\}$ and $T$ is the length of the racing lap. The learned Gaussian process regressors are of the following form:
\begin{equation}
\label{eq:13}
    e_i = \mathscr{GP}(X,Y,\psi,v,\delta,D),\ \ i\in\{X,Y,\psi,v\}.
\end{equation}
where $i$ corresponds to the faulty states defined as $e_Y\thicksim\mathcal{N}(\mu_Y,\sigma_Y)$, $e_i\thicksim\mathcal{N}(\mu_i,\sigma_i)$ $e_X\thicksim\mathcal{N}(\mu_X,\sigma_X)$, $e_{\psi}\thicksim\mathcal{N}(\mu_{\psi},\sigma_{\psi})$ and $e_v\thicksim\mathcal{N}(\mu_v,\sigma_v)$. Each $\mu_i$ and $\sigma_i$ represent the mean and variance of the corresponding distribution. The learned model is used for MPC predictions and can be obtained by adding the prediction model mismatch to model (\ref{eq:1}):
\begin{equation}
\label{eq:15}
    F_\text{learned}(x_k,u_k) = F_\text{mpc}(x_k,u_k)  + F_\text{error}(x_k,u_k).
\end{equation}
The overall approach with MPC controller, MHE estimator, and $\mathscr{GP}$ learning extension is shown in Fig \ref{fig3}.

\section{Results and Discussion}

\subsection{Estimation and \texorpdfstring{$\mathscr{GP}$}{GP} learning results}
\begin{figure}[htb]
\centering
\includegraphics[width=7.75cm,height=5.5cm]{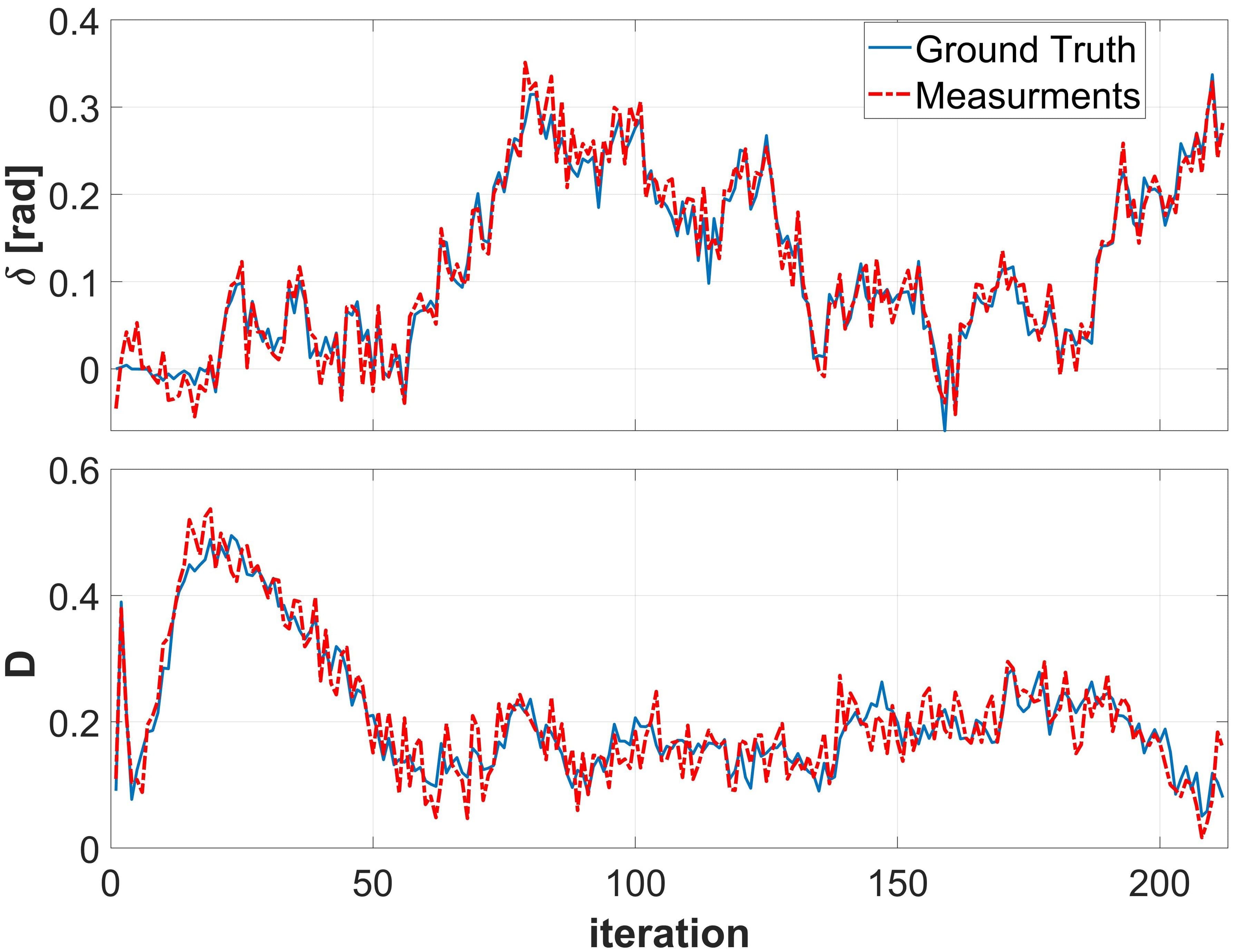}
\caption{Control signals (measurement vs true value).}
\label{fig4}
\end{figure}

\begin{figure}[htb]
\centering
\includegraphics[width=7.75cm,height=5.5cm]{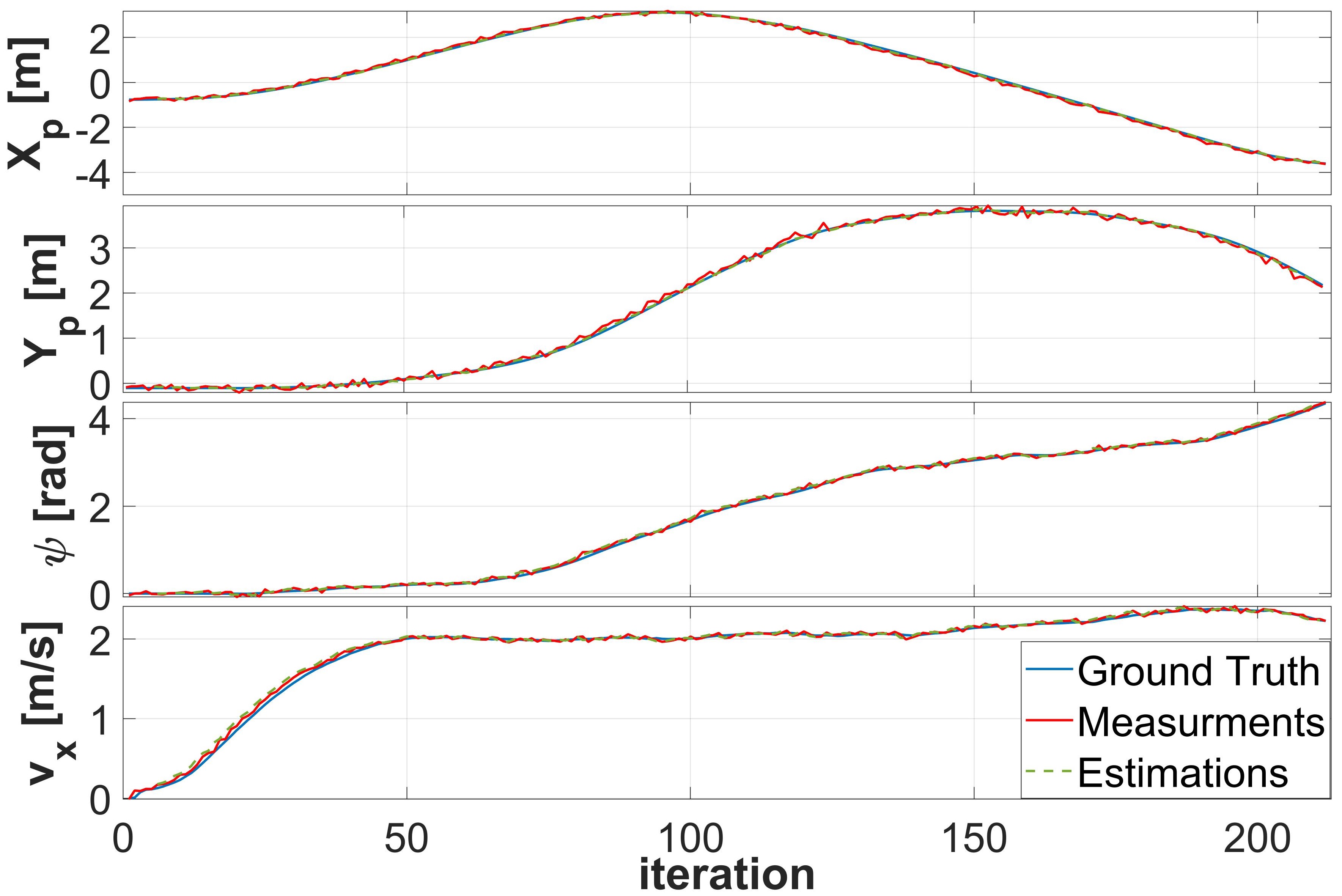}
\caption{Vehicle states estimation with MHE.}
\label{fig5}
\end{figure}

The three optimization problems; MPC for control and planning and MHE for state estimation are coded in \texttt{Matlab} using \texttt{CASADI} framework and solved with \texttt{IPOPT}. Both the MHE estimator and the NMPC controller run at an average frequency of $130Hz$ on a laptop with ryzen $7$ $5800H$ CPU and $32GB$ of DRAM, note that this execution time is more than enough for real-time implementation. The NMPC uses state estimates obtained by the MHE to compute the necessary control actions for tracking the planned racing trajectory. The MHE and NMPC weighting matrices are found by iterative tuning until the desired performance is obtained, they are given along with the measurement noise covariance matrices as follows:
\begin{equation}
    \begin{split}
        Q &= \text{diag}(0.015,0.015,0,0),\\
        R &= \text{diag}(0.003,0.0025),\\
        \sigma_y &= \text{diag}(0.05,0.05,0.035,0.1),\\
        \sigma_u &= \text{diag}(0.2,0.035).  
    \end{split}
\end{equation}
The parameters used for implementation are all properly listed in table \ref{tab2}. Observing Fig. \ref{fig4}, one can notice that despite the low-order model implemented for MHE, the latter was able to estimate the true noiseless states from the noisy states and control measurements (see Fig. \ref{fig5}). The $\mathscr{GP}$ regressors for the faulty states are programmed using sci-kit learn with a combination of radial basis function and constant kernels. The collected dataset is split into 85$\%$ for training and 15$\%$ for testing . The training for the whole lap takes about 15 to 20 seconds and the performance of the $\mathscr{GP}$ regressors is evaluated by the coefficient of determination $R_2$.

\begin{figure}[!h]
\centering
\includegraphics[width=7.75cm,height=11cm]{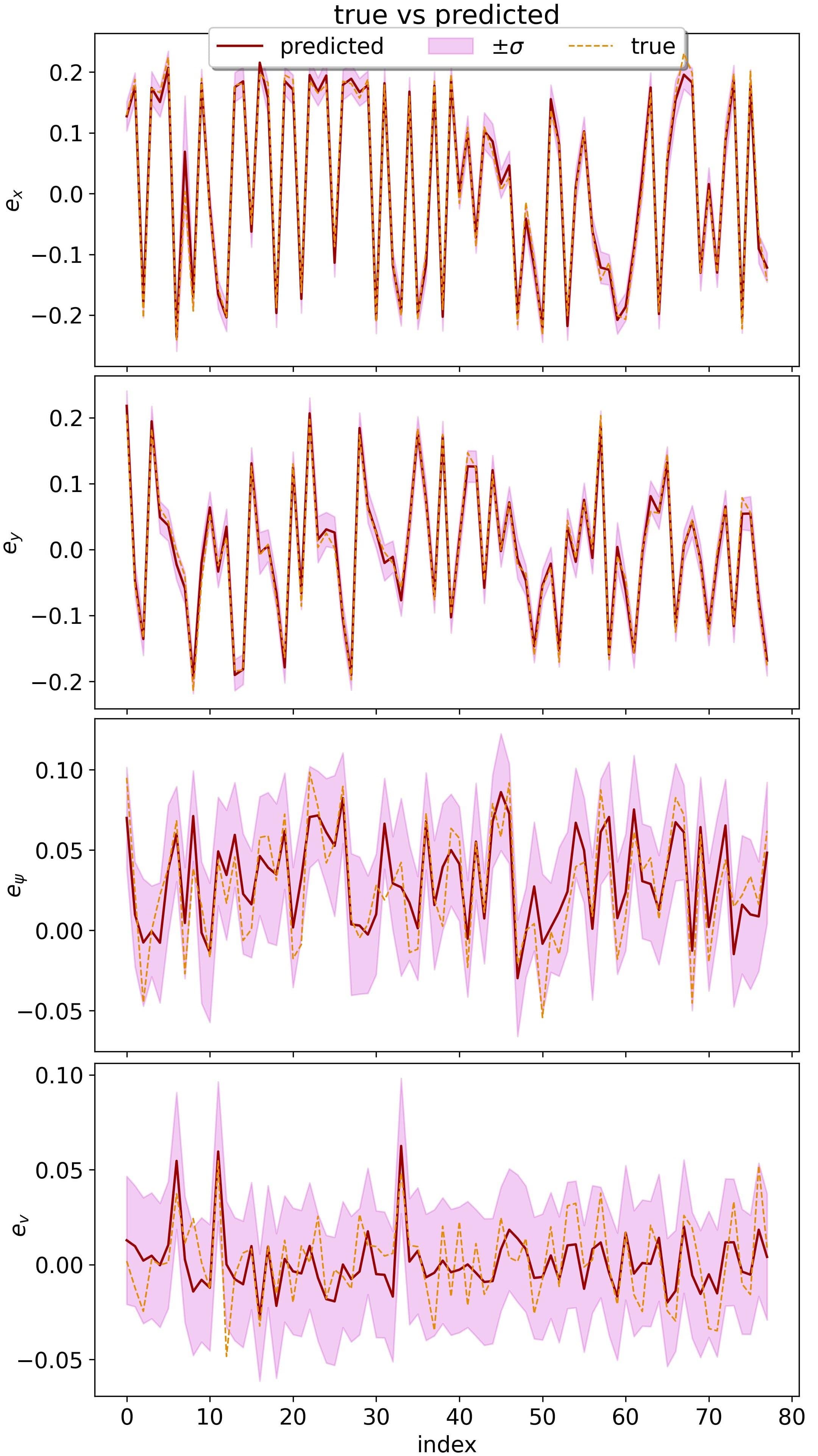}
\caption{\texorpdfstring{$\mathscr{GP}$}{GP} predictions, true values, and confidence intervals for errors of the states (X,\ Y,\ $\psi$ and V).}
\label{fig6}
\end{figure}

Fig. \ref{fig6} shows the $\mathscr{GP}$ predictions and 95$\%$ confidence intervals for the four model mismatch states on the test data. As can be seen, the $\mathscr{GP}$ regressor does a great job in predicting the model mismatch, the highest uncertainties in predictions generally correspond to the track's high-speed corners.

\begin{table}[htb]
\caption{Implementation parameters.} 
\label{tab2}
\centering
\begin{tabular}{c c c c  } 
\hline
Parameter & Value & Parameter & Value  \\ [0.5ex] 
 
$N_e$ & $6$  & $N_p$ & $65$  \\

$V$ & $\sigma_y^{-1}$  & $q$ & $0.025$  \\

$W$ & $\sigma_u^{-1}$  & $r$ & $1.25$ \\

$\delta_{max}$ & $\frac{\pi}{6}$  & $N$ & $16$ \\

$\delta_{min}$ & $-\frac{\pi}{6}$  & $T_s$ & $33 ms$  \\

$D_{max}$ & $1$  & $D_{min}$ & $-1$  \\[1ex] 

\end{tabular}
\end{table}

\subsection{Planning and control results}

\begin{figure}[!t]
\centering
\includegraphics[width=7.75cm,height=5.5cm]{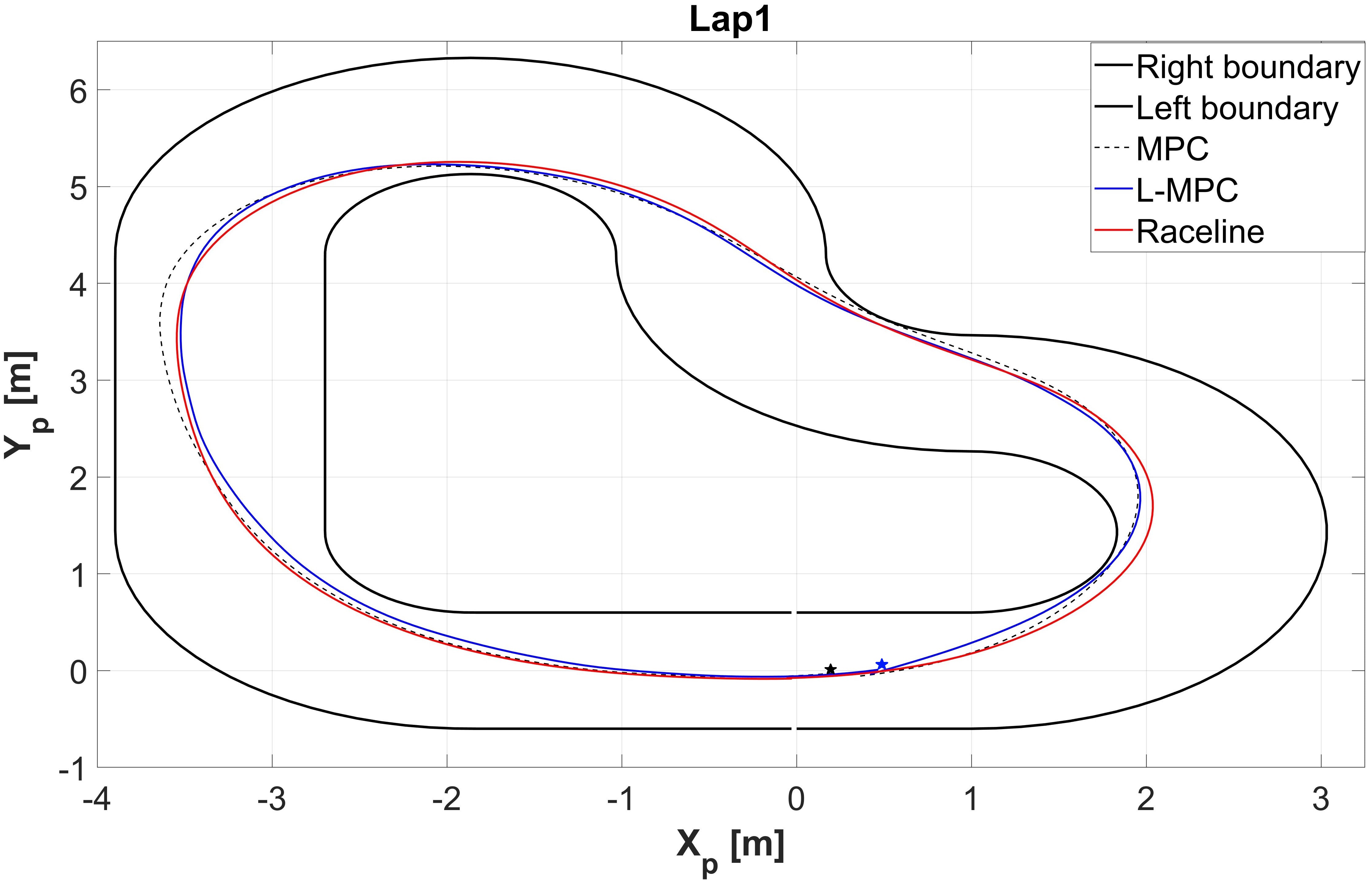}
\caption{Race line tracking for lap 1.}
\label{fig7}
\end{figure}

\begin{figure}[!t]
\centering
\includegraphics[width=7.75cm,height=5.5cm]{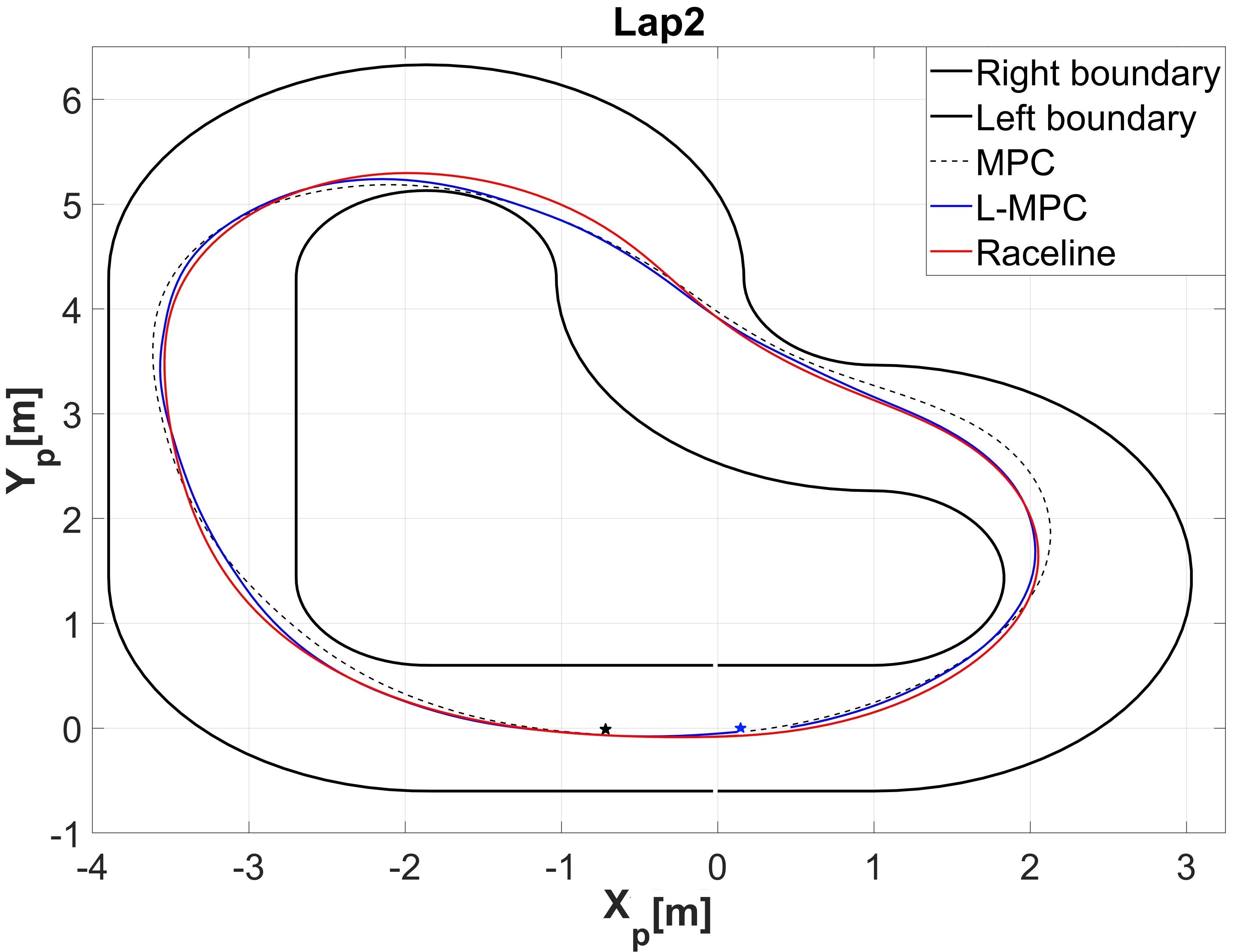}
\caption{Race line tracking for lap 2.}
\label{fig8}
\end{figure}

The designed controller is tested on two race tracks, the first has an L-shape, while the second one is of an oval form. The results of the planning and control algorithms for both L-shaped and oval tracks in two consecutive laps are shown in Figs. \ref{fig7} and \ref{fig8}, and Figs. \ref{fig16} and \ref{fig17}, respectively. The planned race line is represented in red color, and the proposed controller (L-MPC: blue color) is compared to the standard MPC of equation (\ref{eq:2}) (MPC: dotted line). It can be clearly seen that the designed controller tracks well the planned race line while exhibiting acceptable cornering capabilities. In fact, L-MPC has higher tracking precision, where the RMSE error results in 0.091m compared to 0.121m for standard MPC. Moreover, L-MPC performs better cornering compared to regular MPC, which is mainly due to the online model correction with the learning extension allowing for more aggressive driving. Furthermore, the star sign in the figures shows that L-MPC travels a longer distance compared to MPC, meaning that the learning extension allows the controller to reach higher speeds. It is worth noting that tracking precision is further improved in the second lap.

\begin{figure}[!ht]
\centering
\includegraphics[width=7.75cm,height=5.5cm]{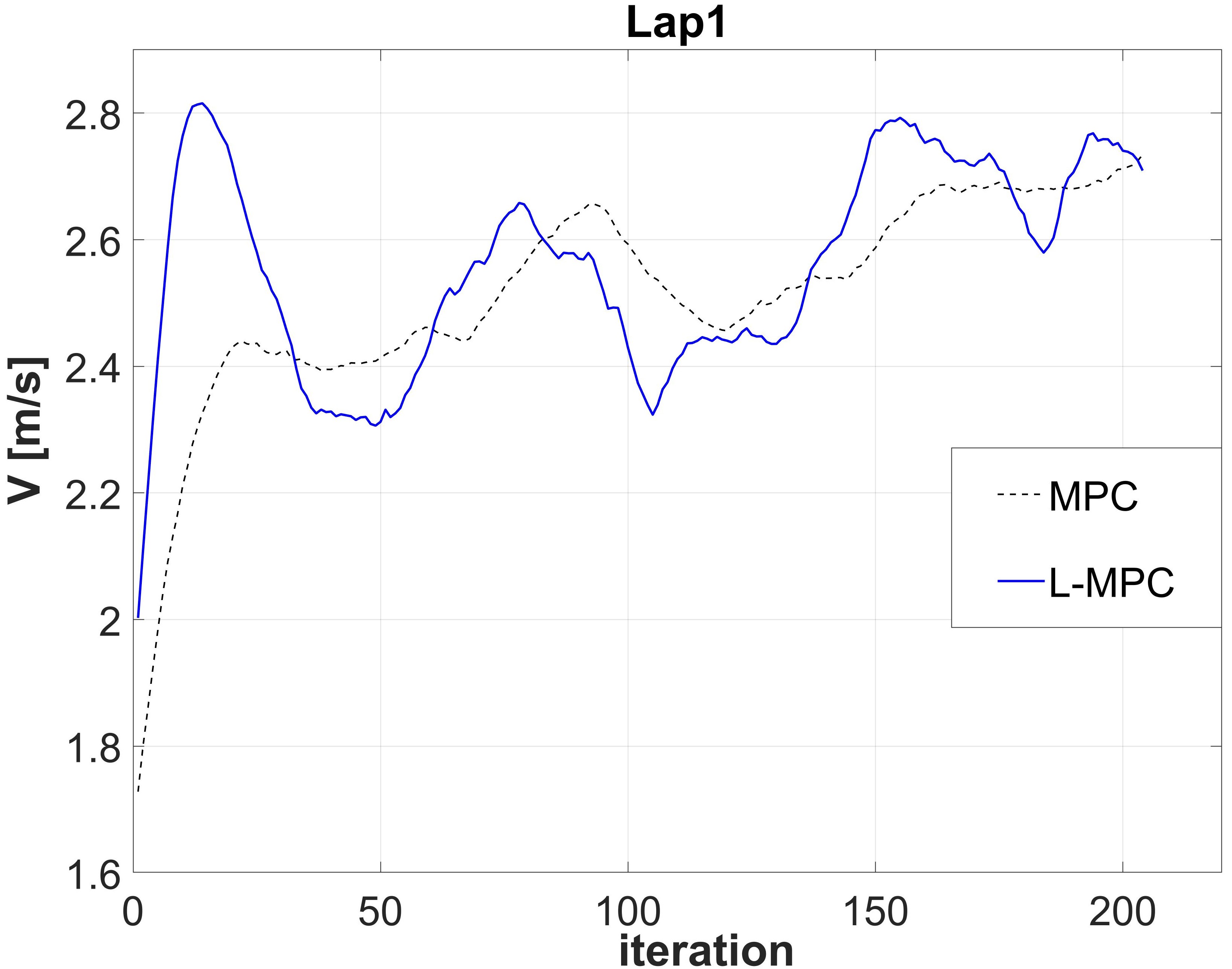}
\caption{Vehicle speed for lap 1.}
\label{fig9}
\end{figure}

\begin{figure}[!ht]
\centering
\includegraphics[width=7.75cm,height=5.5cm]{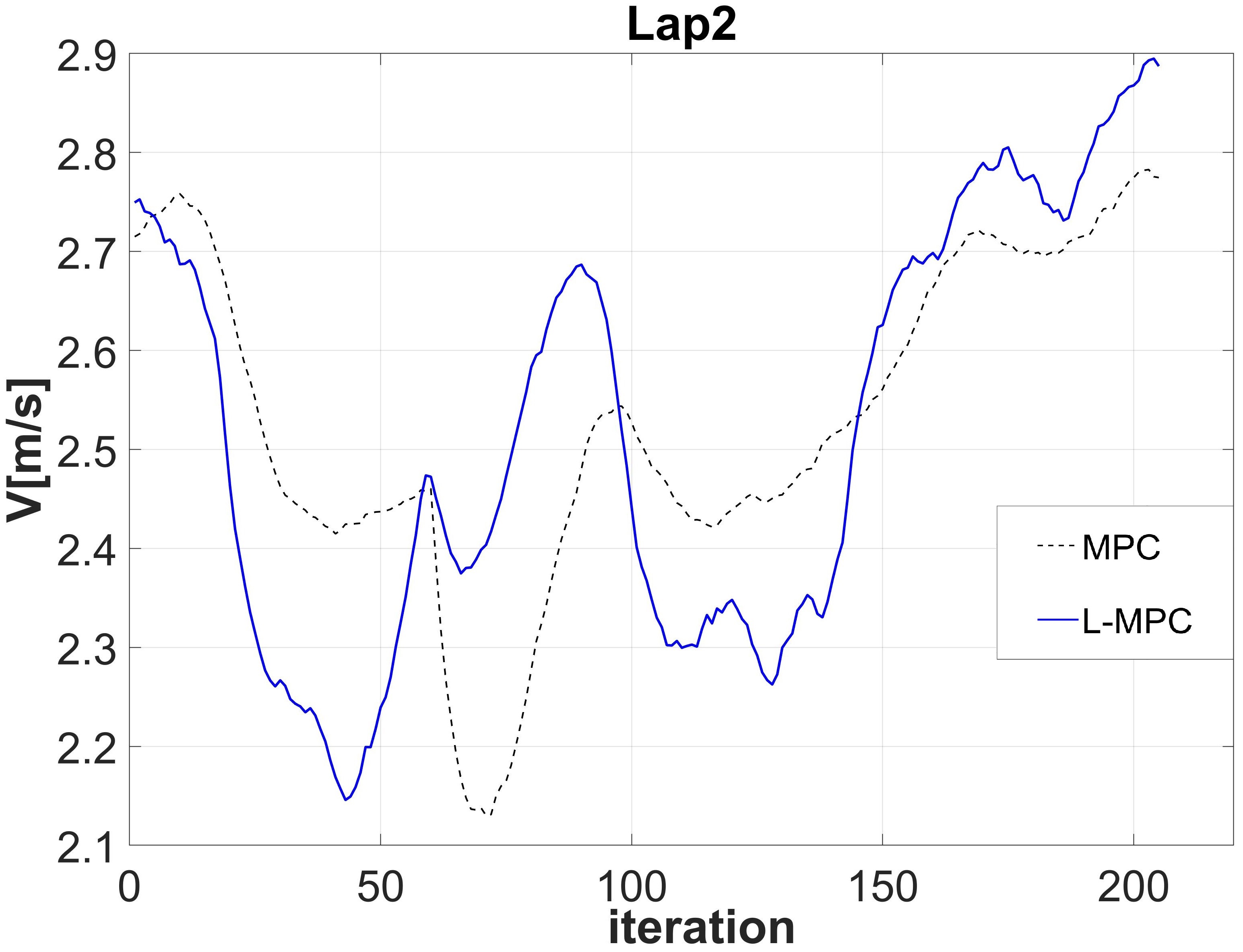}
\caption{Vehicle speed for lap 2.}
\label{fig10}
\end{figure}

\begin{figure}[!ht]
\centering
\includegraphics[width=7.75cm,height=5.5cm]{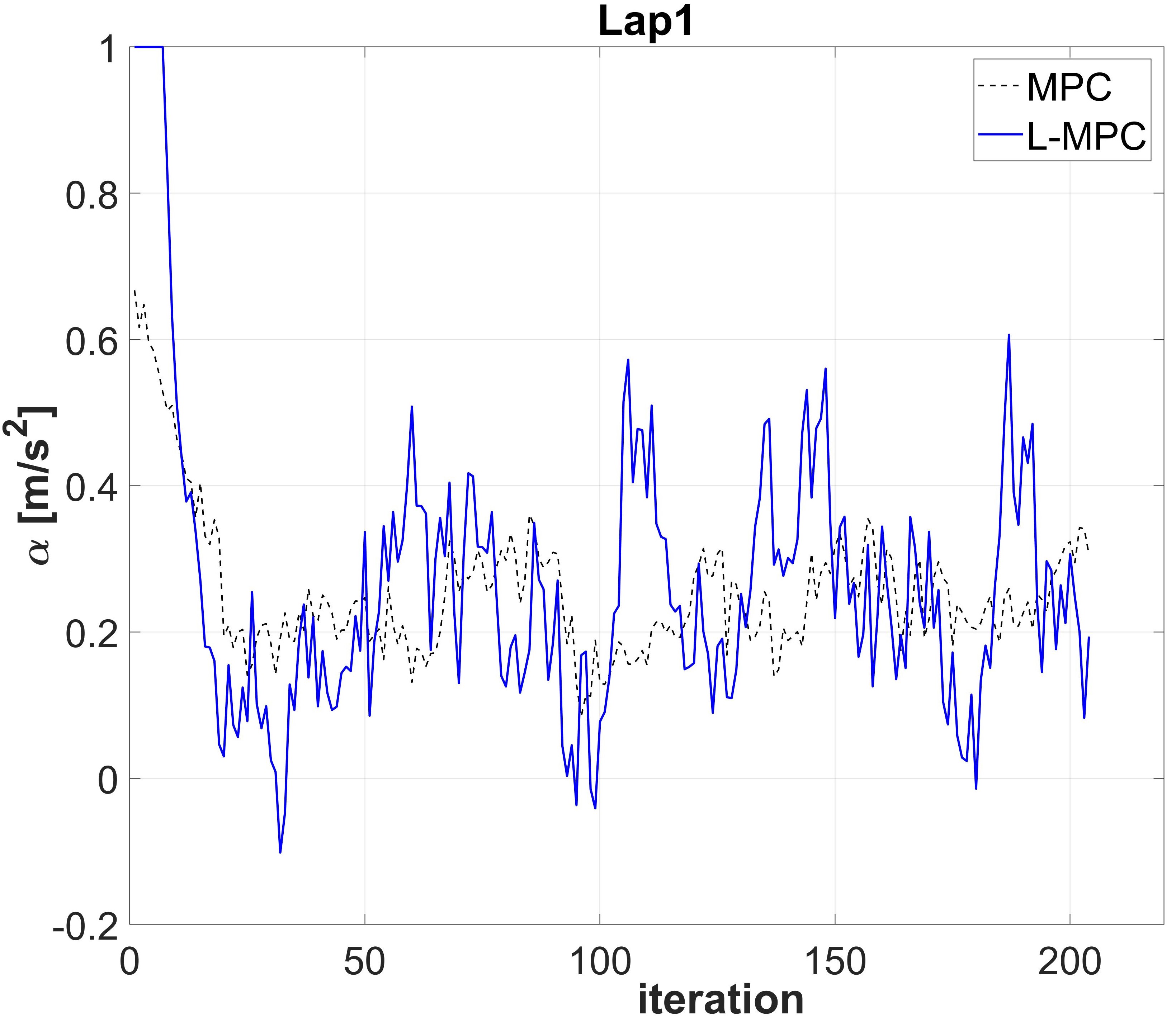}
\caption{Acceleration control for lap 1.}
\label{fig11}
\end{figure}

\begin{figure}[!ht]
\centering
\includegraphics[width=7.75cm,height=5.5cm]{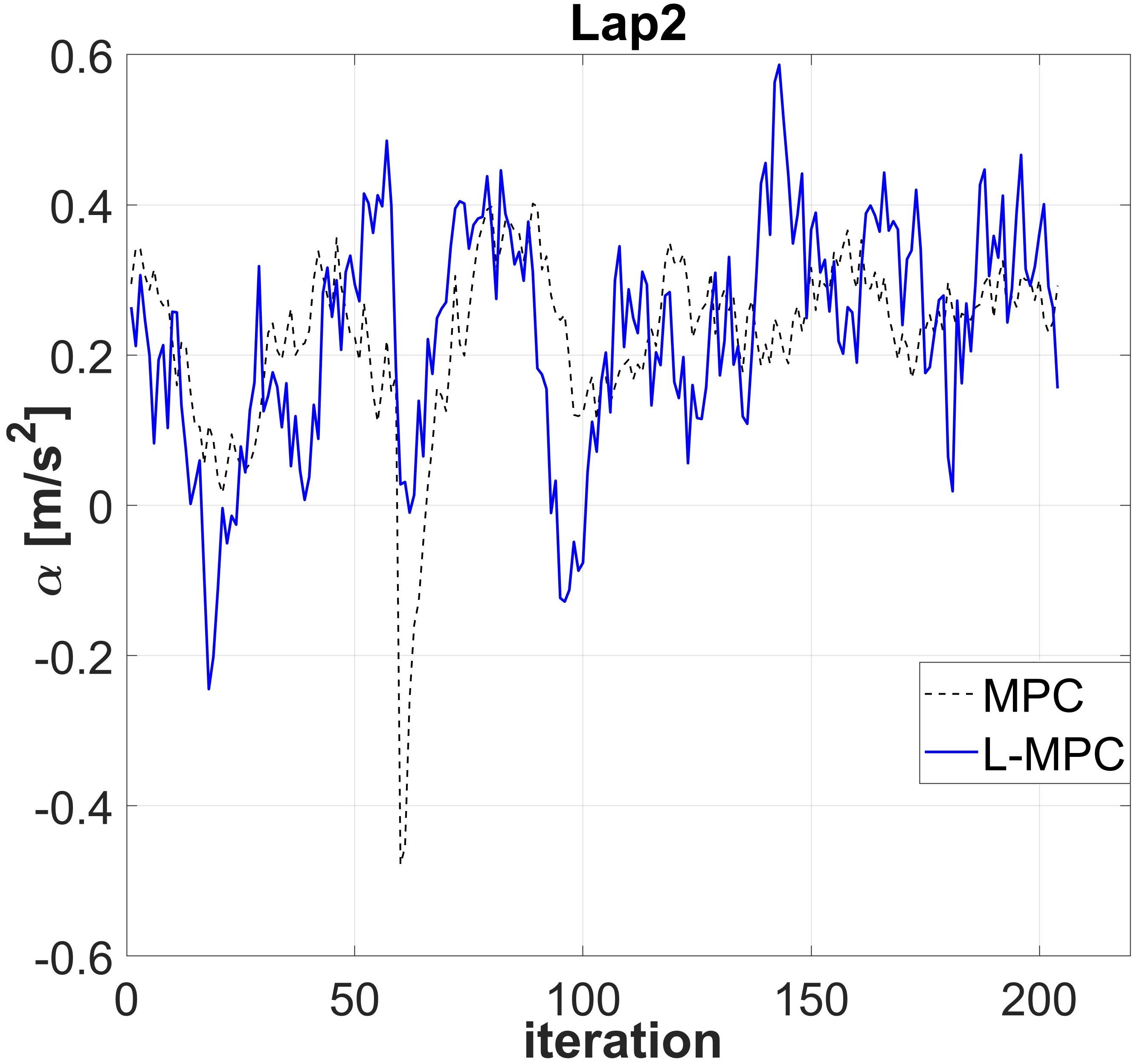}
\caption{Acceleration control for lap 2.}
\label{fig12}
\end{figure}

\begin{figure}[!ht]
\centering
\includegraphics[width=7.75cm,height=5.5cm]{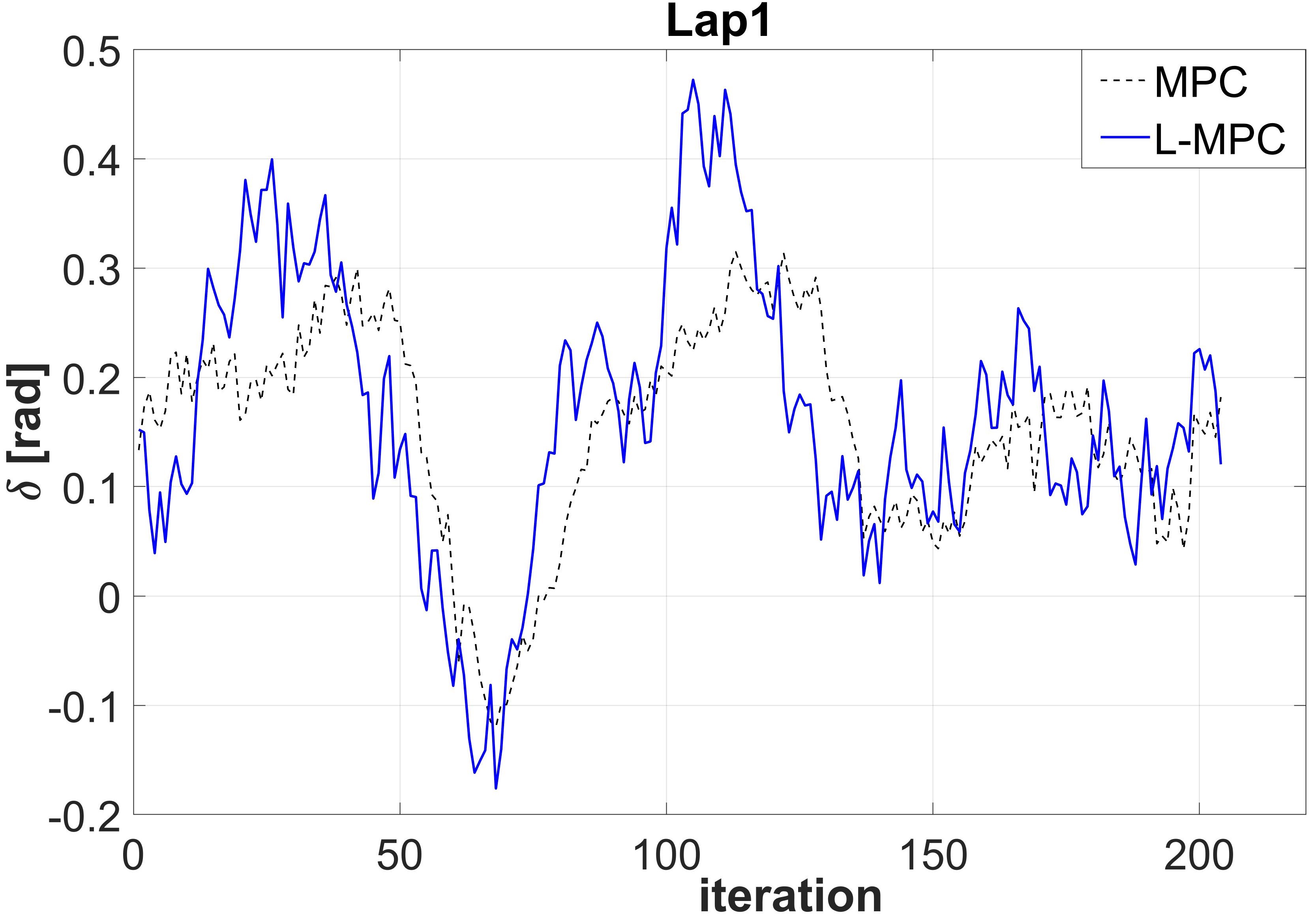}
\caption{Steering control for lap 2.}
\label{fig13}
\end{figure}

\begin{figure}[!ht]
\centering
\includegraphics[width=7.75cm,height=5.5cm]{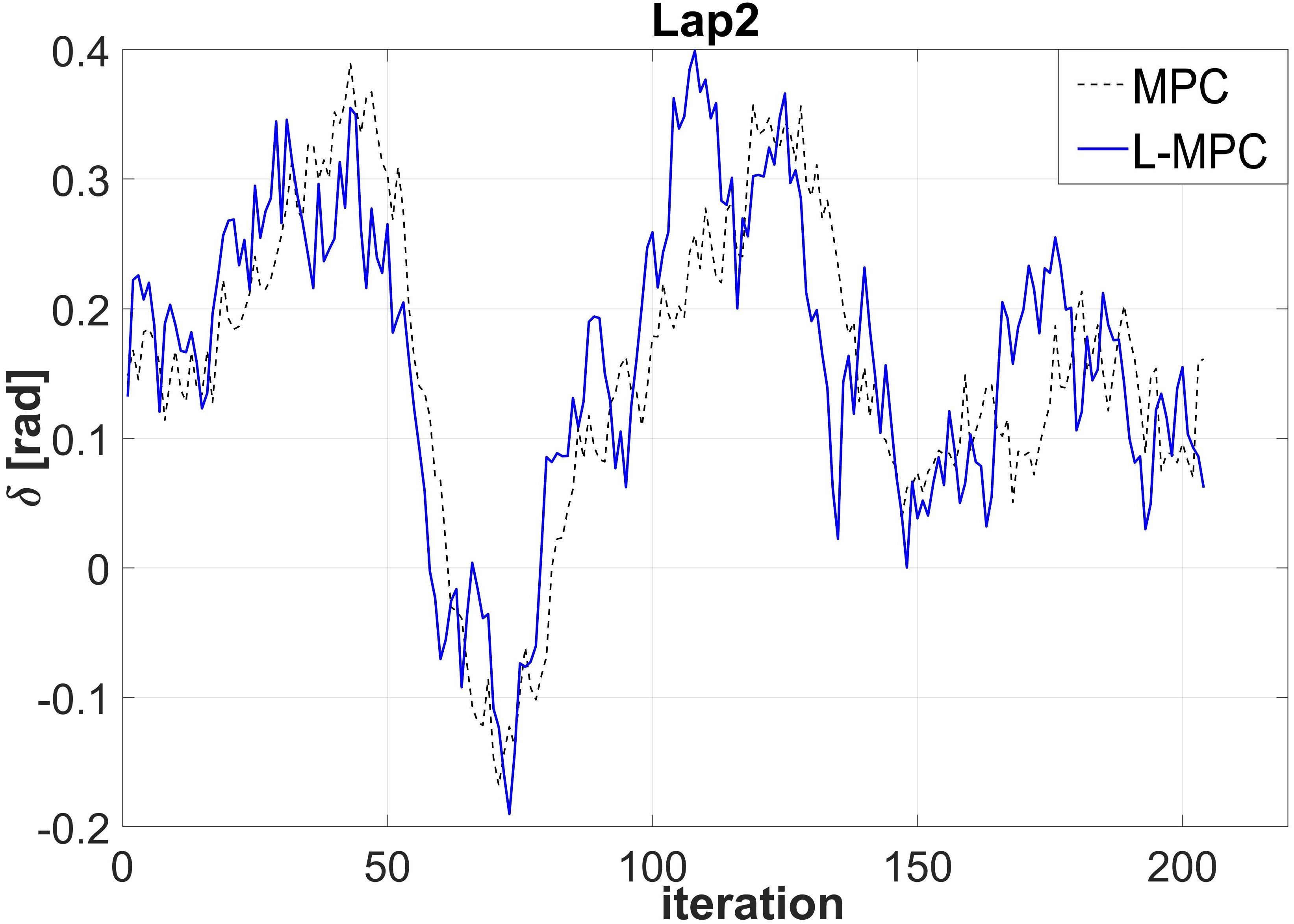}
\caption{Steering control for lap 2.}
\label{fig14}
\end{figure}

Figs. \ref{fig9} and \ref{fig10} show the speed of the vehicle for both laps. Generally, L-MPC is able to reach higher speeds (up to 2.9 m/s), even though regular MPC is also fast enough (up to 2.78 m/s). The average speed for the whole track is evaluated at 2.51 m/s for L-MPC compared to 2.40m/s for MPC. This can be further observed in Figs. \ref{fig11} and \ref{fig12}, where L-MPC exhibits relatively higher accelerations/decelerations and hence more aggressive racing. Figs. \ref{fig13} and \ref{fig14} show the corresponding steering controls of MPC and L-MPC for both laps. 

\begin{figure}[!ht]
\centering
\includegraphics[width=7.75cm,height=5.5cm]{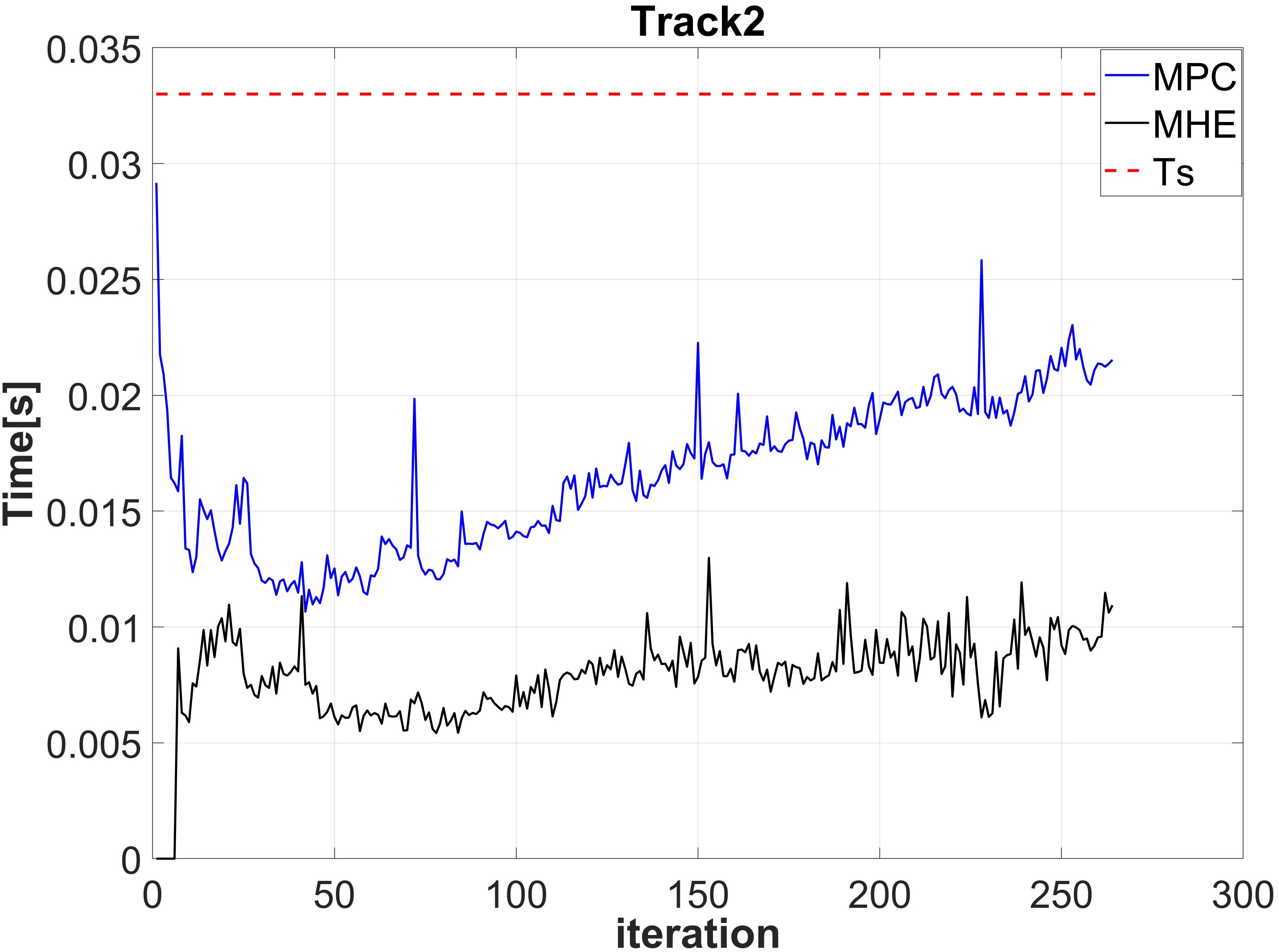}
\caption{Computation time for MPC/MHE.}
\label{fig15}
\end{figure}

The computation time for both MPC and MHE is relatively low, which allows for real-time application. This is very important since autonomous racing systems are known for fast dynamics. More importantly, MHE runs faster than MPC with an average computation time of 9.5ms compared to 20.9ms for MPC. This ensures that the estimates are readily available to the controller and no lag is to be expected. This is clearly illustrated in Fig. \ref{fig15} which shows the required computation time during the first lap of the second track.
The proposed controller has been evaluated on an oval track as well (see Figs. \ref{fig17} and \ref{fig18}), where similar racing performance is observed. Mainly, L-MPC ensures more cornering capabilities with higher tracking precision. The respective tracking errors on the lateral and longitudinal positions are compared in Figs. \ref{fig18} and \ref{fig19} for both laps. The corresponding lateral RMSE (Y position) for the whole track results in 0.044m for L-MPC against 0.082m for standard MPC. Similarly, for the longitudinal errors (X position), L-MPC scored an RMSE of 0.1m against 0.13m for standard MPC. It is worth mentioning that the $\mathscr{GP}$ corrections yielded better improvements in tracking performance for this case compared to the first track, which can be attributed to the challenging nature of the L-shaped track. 

\begin{figure}[!ht]
\centering
\includegraphics[width=7.75cm,height=5.5cm]{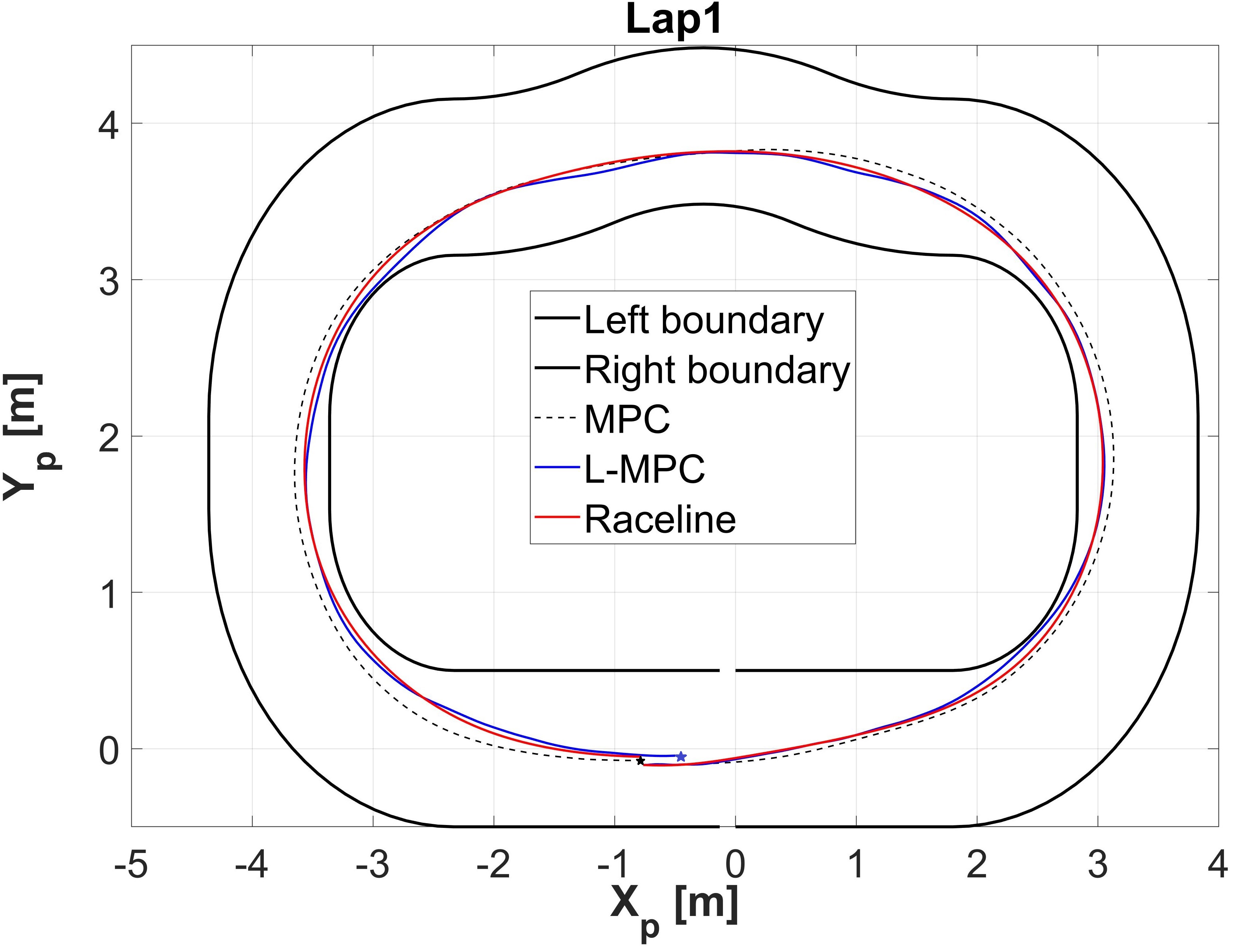}
\caption{Race line tracking for lap 1.}
\label{fig16}
\end{figure}
\begin{figure}[!h]
\centering
\includegraphics[width=7.75cm,height=5.5cm]{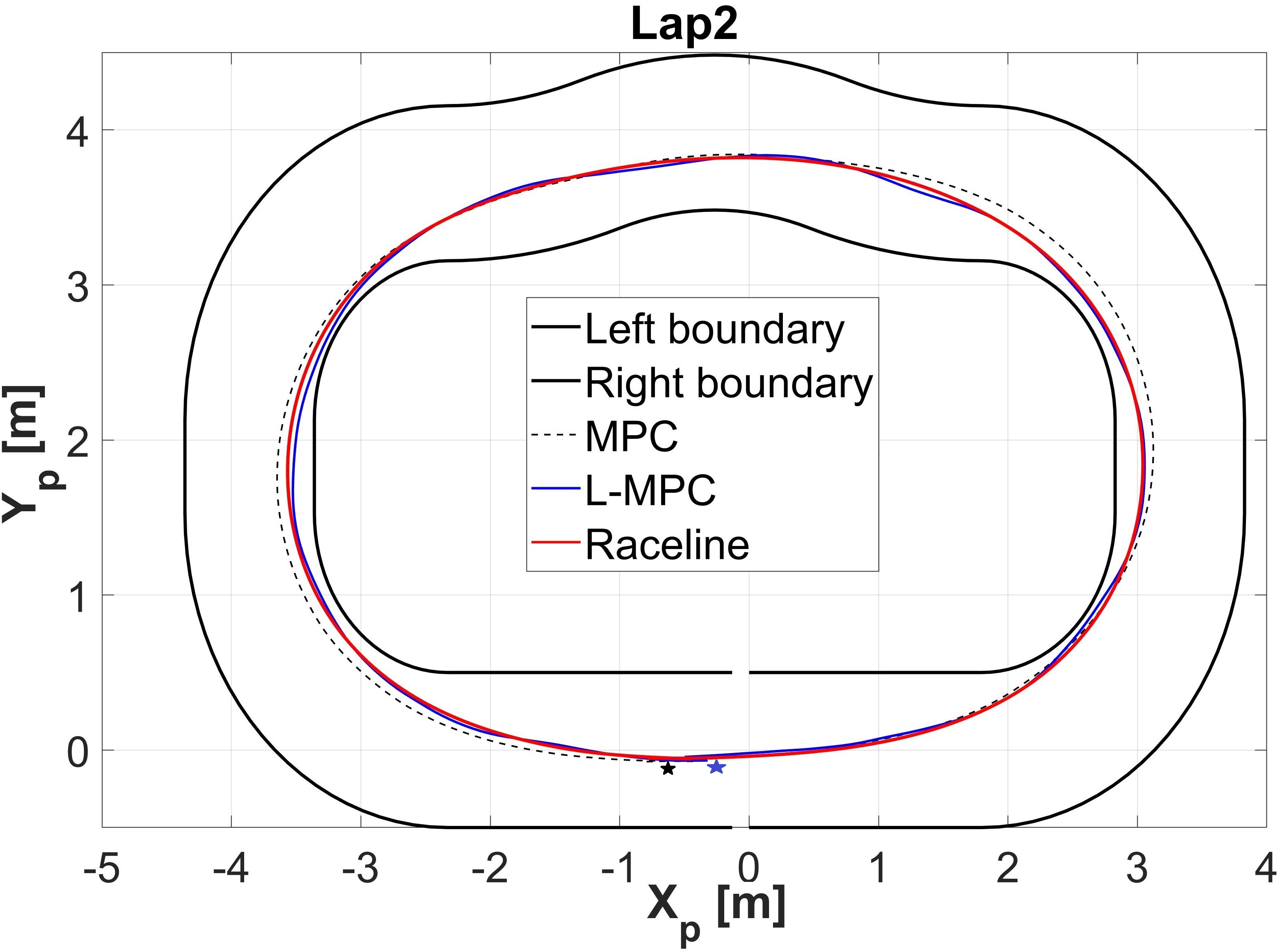}
\caption{Race line tracking for lap 2.}
\label{fig17}
\end{figure}
\begin{figure}[!h]
\centering
\includegraphics[width=7.75cm,height=5.5cm]{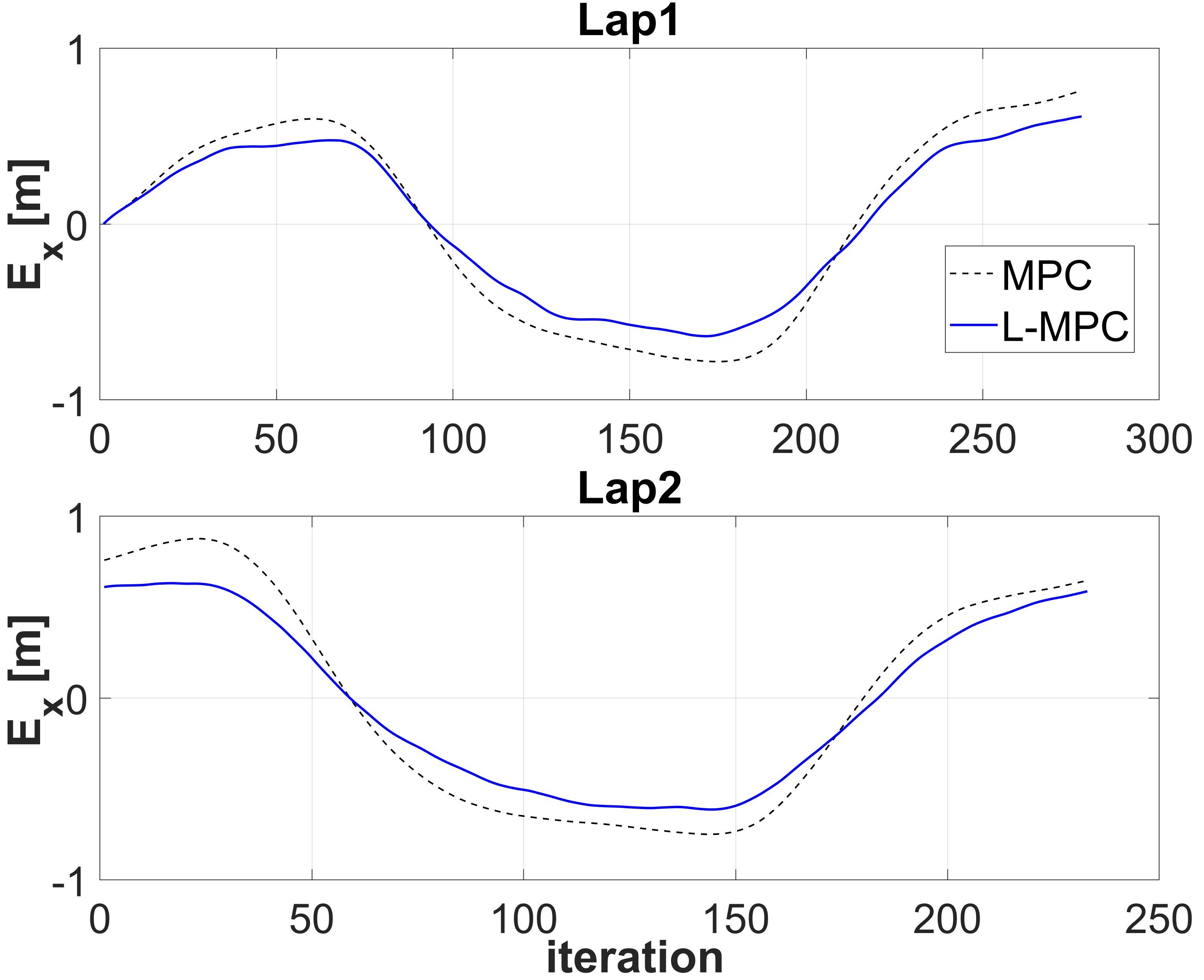}
\caption{Longitudinal tracking error for track 2.}
\label{fig18}
\end{figure}
\begin{figure}[!h]
\centering
\includegraphics[width=7.75cm,height=5.5cm]{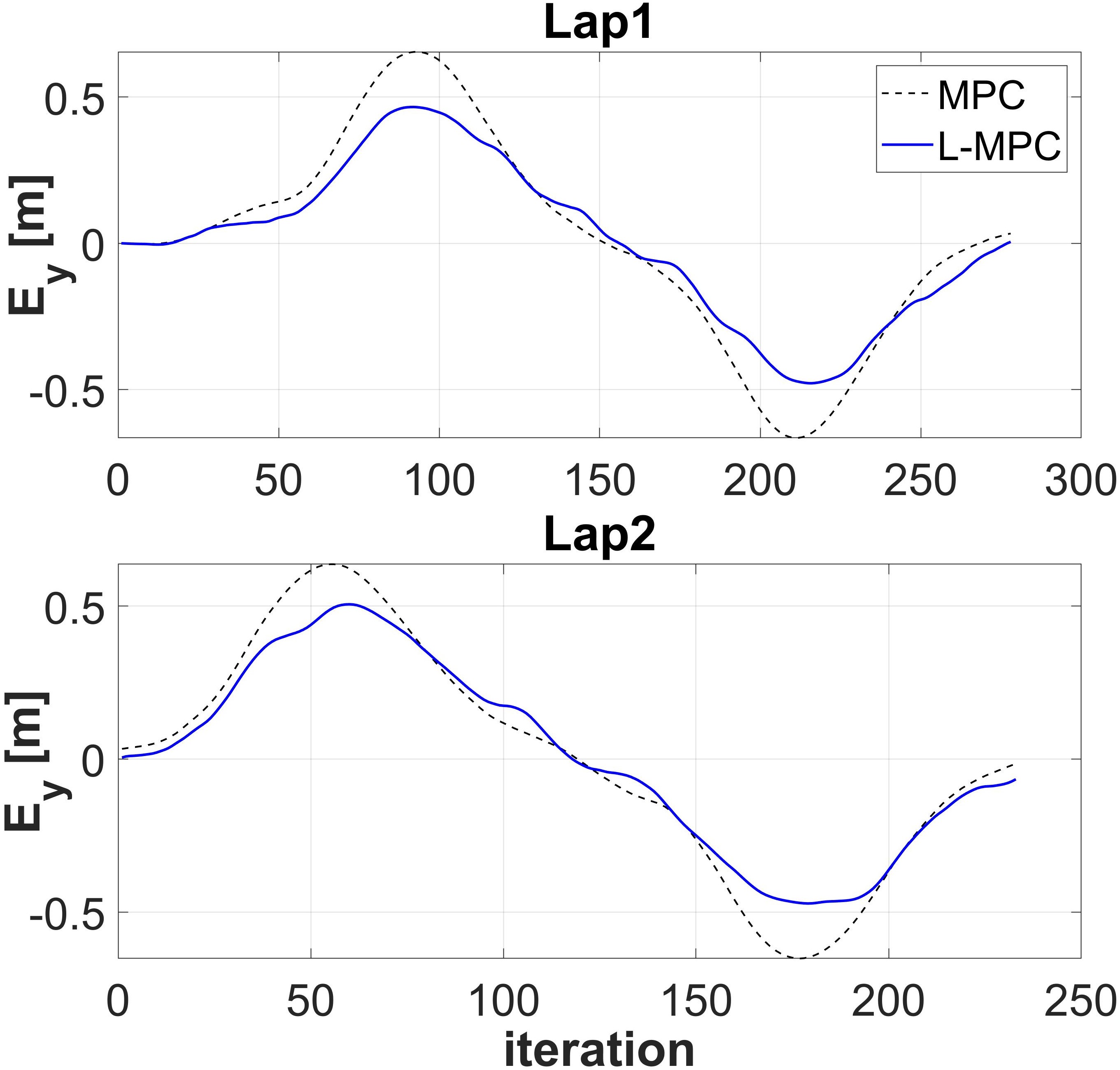}
\caption{Lateral tracking error for track 2.}
\label{fig19}
\end{figure}

\section{Conclusions and Future Works}
In this paper, the problem of autonomous racing has been addressed by developing an NMPC controller that handles both lateral and longitudinal vehicle dynamics. The proposed controller has been coupled with an MHE for state estimation and further improved with a learning extension based on Gaussian process regression to enhance its performance. Furthermore, the optimal race line of the race track has been planned offline using a nonlinear model predictive planner that accounts for vehicle handling limits and circuit boundaries. The proposed controller, estimator, and planner have been tested on two challenging race tracks, and the results proved that the proposed control strategy has superior performance and ensures high-speed racing capabilities. Future works shall address the experimental validation of the proposed approach, and investigate the race line planning approach to be further complemented with obstacle avoidance.

\section*{Acknowledgment}
The authors would like to thank the University of Paris Saclay for the financial support provided to conduct this research.


\bibliographystyle{apalike} 

\bibliography{PhD}

\end{document}